\documentclass[final, a4paper, 11pt]{article}

\usepackage[scale=0.768]{geometry}
\usepackage[english]{babel}

\usepackage[mathscr]{euscript}
\usepackage{color}
\usepackage{fancyhdr}
\usepackage{bm}
\usepackage{hyperref}
\usepackage{cite}
\usepackage{showkeys}
\usepackage{subfigure}
\usepackage{float}

\usepackage{amsfonts}
\usepackage{amsmath}
\usepackage{amssymb}
\usepackage{amscd}
\usepackage{amsthm}

\usepackage{accents}
\usepackage{graphicx}
\usepackage{array}

\usepackage{algorithm}
\usepackage{algorithmicx}
\usepackage{algpseudocode}

\newtheorem{theorem}{Theorem}

\newtheorem{remark}[theorem]{Remark}

\newcommand*{\N}{\ensuremath{\mathbb{N}}}
\newcommand*{\Z}{\ensuremath{\mathbb{Z}}}
\newcommand*{\R}{\ensuremath{\mathbb{R}}}

\renewcommand{\i}{\mathrm{i}}
\renewcommand{\phi}{\varphi}
\renewcommand{\rho}{{\varrho}}
\renewcommand{\epsilon}{{\varepsilon}}

\renewcommand{\d}[1]{\,\mathrm{d}#1 \,}

\newcommand{\ol}[1]{\overline{#1}}
\newcommand{\J}{\mathcal{J}} 
\newcommand{\0}{{0}} 

\renewcommand{\Re}{\mathrm{Re}\,}

\newcommand{\grad}{\nabla}

\newcommand{\W}{{W^*}} 
\newcommand{\Wast}{{W_{\hspace*{-1pt}{\Lambda}^\ast}}} 

\newlength{\dhatheight}

\usepackage{color}

\definecolor{xl}{rgb}{0,0,0}
\definecolor{lxl}{rgb}{0,0,0}

\usepackage{amsmath}
\usepackage{amsfonts}
\usepackage{amssymb}
\begin{document}

\sloppy

\title{Near-field imaging of locally perturbed periodic surfaces}
\author{Xiaoli Liu\thanks{INRIA Saclay Ile de France / CMAP Ecole Polytechnique, Palaiseau, France
;\texttt{xiaoli.liu@inria.fr}} \and 
Ruming Zhang\thanks{Institute of Appied and Numerical mathematics, Karlsruhe Institute of Technology, Karlsruhe, Germany;
\texttt{ruming.zhang@kit.edu}; corresponding author. }
\thanks{The work of the second author was supported by Deutsche Forschungsgemeinschaft (DFG) through CRC 1173.
}}
\date{}
\maketitle

\begin{abstract}
This paper concerns the inverse scattering problem to reconstruct a locally perturbed periodic surface. Different from scattering problems with quasi-periodic incident fields and periodic surfaces, the scattered fields are no longer quasi-periodic. Thus the classical method for quasi-periodic scattering problems no longer works. In this paper, we apply a Floquet-Bloch transform based numerical method to reconstruct both the unknown periodic part and the unknown local perturbation from the near-field data. 

By transforming the original scattering problem into one defined in an infinite rectangle, the information of the surface is included in the coefficients. The numerical scheme contains two steps. The first step is to obtain an initial guess, i.e., the locations of both the periodic surfaces and the local perturbations, from a sampling method. The second step is to reconstruct the surface. As is proved in this paper, for some incident fields, the corresponding scattered fields  carry little information of the perturbation. In this case, we use this scattered field to reconstruct the periodic surface. Then we could apply the data that carries more information of the perturbation to reconstruct the local perturbation. The Newton-CG method is applied to solve the associated  optimization problems.  Numerical examples are given at the end of this paper to show the efficiency of the numerical method.
\end{abstract}

\section{Introduction}

In this paper, we introduce the numerical method of the inverse scattering problem from a locally perturbed periodic surface. Both the periodic part and the local perturbation  of the surface are unknown. {The} aim of the inveres problem is to reconstruct both of them from the {near-filed measurement} data. 

 {Since} the periodic surface is perturbed, the classical framework for the quasi-periodic scattering problems (i.e., quasi-periodic incident fields with periodic domains) no longer works. An efficient way to solve these problems is to apply the Floquet-Bloch transform. With the help of this Fourier-like transform, the original problem, which is defined in a 2D unbounded domain, is written into a new one defined in a 3D bounded domain. This method has been applied to perturbed periodic structures in \cite{Coatl2012} and waveguide problems in \cite{Hadda2016}.  For scattering problems with non-periodic incident fields and periodic surfaces, we refer to \cite{Lechl2016a,Lechl2016b}. For problems with locally perturbed periodic surfaces, see \cite{Lechl2016,Lechl2017}. In the paper \cite{Zhang2017e}, a high order numerical method has been proposed based on the Floquet-Bloch transform {and} this method is used in this paper to produce the measured data. For a fast imaging method to reconstruct the local perturbations in periodic media with the help of the Bloch transform, we refer to \cite{Cakon2018}.

 The work in this paper is an extension {to} the joint work of the second author with {\em Prof. Armin Lechleiter} in \cite{Zhang2017c}, where the periodic surface is assumed to be already known. A numerical method has been proposed to find out both the location and the shape of the local perturbation. The sampling method introduced by Ito, Jin and Zou (see \cite{Ito2012a}) was extended to  find out the location, and a Newton-CG method was applied to reconstruct the shape. However, the setting in this paper is more difficult, i.e., the periodic surface is no longer known. Thus we have to find out the location without a known periodic surface, and also reconstruct both the periodic surface and the perturbation. {In this case}, the sampling method in \cite{Zhang2017c} {does not work any more}, which makes the problem much more challenging.

In this paper, we {develop} a numerical method for the inverse problem. The first task is to find out the locations of both the perturbation and the periodic surface. {Since} the {previous sampling} method does not work, we apply the rough surface reconstruction algorithm introduced in \cite{Liu2018} to obtain the locations, when  {the perturbation} is assumed to be existing in a relatively large domain. Then we apply the Newton's method to reconstruct the shapes of both the periodic surface and the local perturbation. 
The reconstruction contains two steps. The first step is to reconstruct the periodic surface. From an estimate of the difference between field with and without perturbation, for certain incident fields, {the scattered field with perturbation could be a good approximation of the one without perturbation}. Thus in this case, the measured scattered field could be adopted to reconstruct the periodic surface. Base on the {\color{lxl}former} approximation of the periodic surface, we {\color{lxl}apply} the method in \cite{Zhang2017c} to find out the approximation of the local perturbation.

The rest of the paper is organized as follows. In Section 2, we {\color{lxl}recall} the mathematical model of the direct scattering problem and the Floquet-Bloch transform based formulation. In Section 3, the estimation is considered for the difference between scattered fields with and without perturbation. The inverse problem is formulated in Section 4, and the Fr\'{e}chet derivative and its adjoint operator are studied. In Section 5, we conclude the algorithm for inverse problems, including the sampling method for the initial guess and the iterative method for the reconstruction. In Section 6,  we present two numerical results {\color{lxl}obtained from our} algorithm.


\section{Direct Scattering Problem}
\subsection{Mathematical Model}

Given a bounded $2\pi$-periodic function $\zeta$, it defines a periodic surface
\begin{equation*}
 \Gamma:=\left\{(x_1,\zeta(x_1)):\,x_1\in\R\right\}.
\end{equation*}
Let the function $p$ be a compactly supported perturbation. For simplicity, suppose ${\rm supp}(p)\subset (-\pi,\pi)+2\pi J$, where $J\in\Z$ is an integer.
Let $\zeta_p:=\zeta+p$ be the perturbed function {\color{lxl}and} define
\begin{equation*}
 \Gamma_p:=\left\{(x_1,\zeta_p(x_1)):\, x_1\in\R\right\}.
\end{equation*}
Let the domain above $\Gamma$ be $\Omega$ and that above $\Gamma_p$ be $\Omega_p$.

\begin{figure}[H]
\centering
 \includegraphics[width=15cm]{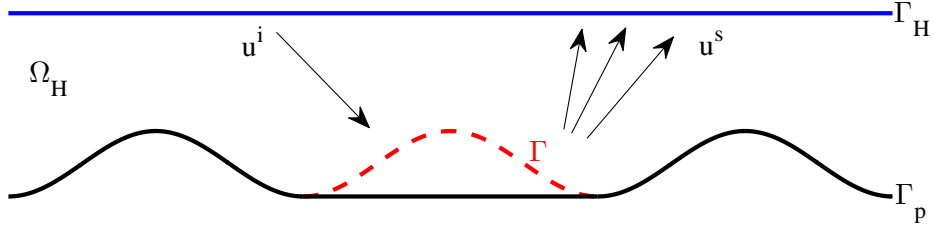}
 \caption{Mathematical model for the scattering problem.}
 \label{fig:model}
\end{figure}

\begin{remark}
For simplicity, from Section 2 to Section 3, we fix $J=0$ for  theoretical arguments.
\end{remark}

In this paper, we assume that the surface $\Gamma_p$ is sound-soft.
Given an incident field $u^i$ that satisfies $\Delta u^i+k^2 u^i=0$ in $\Omega_p$, then it is scattered by $\Gamma_p$ and {\color{lxl}generates} the scattered field {\color{lxl}$u^s$} (or equivalently, the total field $u=u^i+u^s$). For the mathematical model we refer to Figure \ref{fig:model}. First, $u$ satisfies
\begin{equation}\label{eq:sca1}
  \Delta u+k^2 u=0\quad\text{ in }\Omega_p.
\end{equation}
Second, as the surface $\Gamma_p$ is sound-soft, 
\begin{equation}\label{eq:sca2}
 u=0\quad\text{ on }\Gamma_p.
\end{equation}
Moreover, the scattered field $u^s$ is propagating upwards. The Upward Propogation Radiation Condition (UPRC) is typically written as a double layer potential, see \cite{Chand1998a}, and an alternative definition was introduced in \cite{Chand2005,Chand2010}. Let $H$ be a real number that is larger than $\|\zeta\|_\infty$ and $\|\zeta_p\|_\infty$. Then the UPRC is written as
\begin{equation*}
 u^s(x_1,x_2)=\frac{1}{2\pi}\int_\R e^{\i\xi x_1+\i\sqrt{k^2-\xi^2}\,(x_2-H)}\widehat{u}^s(\xi,H)\d\xi,\quad x_2\geq H,
\end{equation*}
where $\widehat{u}^s(\xi,H)$ is the Fourier transform of $u^s(x_1,H)$. Define the Dirichlet-to-Neumann map $T^+$ by
\begin{equation*}
 \left(T^+\phi\right)(x_1)=\frac{\i}{2\pi}\int_\R \sqrt{k^2-\xi^2}\,e^{\i\xi x_1}\widehat{\phi}(\xi)\d\xi,\quad \phi=\frac{1}{2\pi}\int_\R e^{\i \xi x_1}\widehat{\phi}(\xi)\d\xi.
\end{equation*}
Let $ \Gamma_H:=\left\{(x_1,H):\, x_1\in\R\right\},$ then the UPRC is equivalent to 
\begin{equation}\label{eq:sca3}
 \frac{\partial u}{\partial x_2}(x_1,H)=T^+\left[u\big|_{\Gamma_H}\right]+f\,\text{ on }\Gamma_H,\,\text{where}\,\, f=\frac{\partial u^i}{\partial x_2}(x_1,H)-T^+\left(u^i\big|_{\Gamma_H}\right).
\end{equation}

{\color{lxl}Define} the domain $\Omega^p_H:=\R\times(-\infty,H)\cap\Omega_p$ {\color{lxl}and} we consider the problems \eqref{eq:sca1}-\eqref{eq:sca3} in the weighted Sobolev space $H_r^1(\Omega^p_H)$, {\color{lxl}where the} space $H_r^1(\Omega^p_H)$ {\color{lxl}is} defined by
\begin{equation}
 H_r^1(\Omega^p_H):=\left\{\phi\in\mathcal{D}'(\Omega^p_H):\,(1+|x|^2)^{r/2}\phi\in H^1(\Omega^p_H)\right\}.
\end{equation}
{\color{lxl}$\widetilde{H}_r^1(\Omega^p_H)$} is the subspace of $H_r^1(\Omega^p_H)$ such that all the elements vanish on $\Gamma_p$. 
Similarly, we can {\color{lxl}define} the weighted spaces $H^{1/2}_r(\Gamma_H)$ and $H^{-1/2}_r(\Gamma_H)$.  From \cite{Chand2010}, the operator $T^+$ is bounded and continuous from $H_r^{1/2}(\Gamma_H)$ to $H_r^{-1/2}(\Gamma_H)$ for all $|r|<1$. 

The weak formulation of the scattering problem \eqref{eq:sca1}-\eqref{eq:sca3} {\color{lxl}is to find}  $u\in \widetilde{H}_r^1(\Omega^p_H)$ such that
\begin{equation}\label{eq:sca_var}
 \int_{\Omega^p_H}\left[\grad u\cdot\grad\overline{v}-k^2 u\overline{v}\right]\d x-\int_{\Gamma_H}T^+\left(u\big|_{\Gamma_H}\right)\overline{v}\d s=\int_{\Gamma_H}f\overline{v}\d s
\end{equation}
for all $v\in \widetilde{H}^1_r(\Omega^p_H)$ with compact support in $\overline{\Omega^p_H}$. The unique solvability of the variational problem \eqref{eq:sca_var} has been proved in \cite{Chand2010}:.

\begin{theorem}\label{th:uni_sca}
 For $|r|<1$, given any incident field $u^i$ and the function $f$ defined by \eqref{eq:sca3} belongs to the space $ H_r^{-1/2}(\Gamma_H)$, the variational problems \eqref{eq:sca_var} has a unique solution $u\in \widetilde{H}^1_r(\Omega^p_H)$.
\end{theorem}

\begin{remark}
 Although the unique solvability is proved for bounded surfaces, in this paper, the functions $\zeta$ and $\zeta_p$ are assumed to be at least Lipschitz continuous.
\end{remark}

\subsection{Floquet-Bloch transform}

During the numerical process of the inverse {\color{lxl}problem}, a Floquet-Bloch transform based numerical method is applied to solve the direct scattering problems. Thus in this section, we give a brief introduction to this method. Let $h$ and $H$ be two real numbers such that $h<\min\{\zeta,\zeta_p\}<\max\{\zeta,\zeta_p\}<H$, then define $D:=\R\times(h,H)$. Define the periodic cell $W$ and its dual-cell $\W$ by
\begin{equation*}
 W=(-\pi,\pi],\quad \W=(-1/2,1/2].
\end{equation*}
Then let $D^{2\pi}=D\cap W\times\R$, $\Gamma^{2\pi}_h=W\times\{h\}$, $\Gamma_H^{2\pi}=W\times\{H\}$. 
Define the Bloch transform with period $2\pi$ in $D$ by
\begin{equation*}
 \J_D\phi(\alpha,x)=\sum_{j\in\Z}\phi\left(x+\left(
\begin{matrix}
 2\pi j\\0
\end{matrix}
\right)\right)e^{2\i\pi j\alpha}.
\end{equation*}
Define the function space $H_0^r(\W;H^s_\alpha(D^{2\pi}))$ by the closure of $C_0^\infty(\W\times D^{2\pi})$ with the following norm for $r\in\N$:
\begin{equation*}
 \|\phi\|_{H_0^r(\Wast;H^s_\alpha(D^{2\pi}))}=\left[\sum_{\gamma=0}^r\int_\W \|\partial^\gamma_\alpha\psi(\alpha,\cdot)\|^2_{H_\alpha^s(D^{2\pi})}\right]^{1/2}.
\end{equation*}
The definition is extended to all $r\geq 0$ by interpolation between Hilbert space, and to $r\in\R$ by duality arguments. The property of the Bloch transform has been investigated in \cite{Lechl2017}.
\begin{theorem}\label{th:Bloch}
 The Bloch transform is an isomorphism between $H_r^s(D)$ and $H_0^r(\W;H^s_\alpha(D^{2\pi}))$. Further, when $s=r=0$, $\J_D$ is an isometry with the inverse
 \begin{equation}
  (\J_D^{-1}\phi)\left(x+\left(\begin{matrix}
                                     2\pi j\\ 0
                                    \end{matrix}
  \right)\right)=\int_\W \phi(\alpha,x)e^{2\i\pi j\alpha}\d\alpha,\quad x\in D^{2\pi},
 \end{equation}
and the inverse transform equals to the adjoint operator of $\J_D$.
\end{theorem}

Now we apply the Floquet-Bloch transform to the scattering problem \eqref{eq:sca_var}. Following \cite{Lechl2017}, the first task is to transform the original problem, which is defined in the non-periodic domain $\Omega_H^p$, to a periodic domain. In this paper, we choose $D$ as the periodic domain. Let $H_0$ be a real number that lies in the interval $(\min\{\zeta,\zeta_p\},H)$
{\color{lxl}and} define the following two diffeomorphisms for $x\in\Omega^p_{H_0}$:
\begin{align*}
 \Phi_\zeta:\,x\mapsto\left(x_1,x_2+\frac{(x_2-H_0)^3}{(h-H_0)^3}(\zeta(x_1)-h)\right);
 \quad\Phi_{\zeta_p}:\,x\mapsto\left(x_1,x_2+\frac{(x_2-H_0)^3}{(h-H_0)^3}(\zeta_p(x_1)-h)\right).
\end{align*}
Then extend them by the identity operator for $x_2\geq H_0$. {From the assumption that ${\rm supp}(\zeta_p-\zeta)\subset W$, ${\rm supp}(\Phi_\zeta-\Phi_{\zeta_p})\subset D^{2\pi}$. }

 Let $u_D=u\circ\Phi_{\zeta_p}$, {\color{lxl}it} is easily checked that $u_D$ {\color{lxl}satisfies the}  following variational equation:
\begin{equation}\label{eq:var_trans}
 \int_{D}\left[A_{\zeta_p}\grad u_D\cdot\grad\overline{v_D}-k^2 c_{\zeta_p} u_D\overline{v_D}\right]\d x-\int_{\Gamma_H}T^+\left(u_D\big|_{\Gamma_H}\right)\overline{v_D}\d s=\int_{\Gamma_H}f\overline{v_D}\d s,
\end{equation}
for all $v_D=v\circ\Phi_{\zeta_p}\in \widetilde{H}^1(D)$, where
\begin{eqnarray*}
 && A_{\zeta_p}(x)=\left|\det\grad\Phi_{\zeta_p}(x)\right|\left[\left(\grad\Phi_{\zeta_p}(x)\right)^{-1}\left(\left(\grad\Phi_{\zeta_p}(x)\right)^{-1}\right)^\top\right]\in L^\infty(D,\R^{2\times 2});\\
 &&c_{\zeta_p}(x)=\left|\det\grad\Phi_{\zeta_p}(x)\right|\in L^\infty(D).
\end{eqnarray*}
We define the matrix $A_\zeta$ and $c_\zeta$ by $\Phi_\zeta$ in the similar way, i.e.,
\begin{eqnarray*}
 && A_{\zeta}(x)=\left|\det\grad\Phi_{\zeta}(x)\right|\left[\left(\grad\Phi_{\zeta}(x)\right)^{-1}\left(\left(\grad\Phi_{\zeta}(x)\right)^{-1}\right)^\top\right]\in L^\infty(D,\R^{2\times 2});\\
 &&c_{\zeta}(x)=\left|\det\grad\Phi_{\zeta}(x)\right|\in L^\infty(D).
\end{eqnarray*}
As ${\rm supp}(\Phi_\zeta-\Phi_{\zeta_p})\subset D^{2\pi}$, the supports of both $A_{\zeta_p}-A_\zeta$ and $c_{\zeta_p}-c_\zeta$ are subsets of $D^{2\pi}$. Let $w:=\J_D u_D$, then it satisfies
\begin{equation}\label{eq:var_bloch}
  \int_\W a_\alpha(w(\alpha,\cdot),z(\alpha,\cdot))\d\alpha+b(\J_\Omega^{-1}w,\J_\Omega^{-1}z)=\int_\W\int_{\Gamma^{2\pi}_H}F(\alpha,\cdot)\overline{z(\alpha,\cdot)}\d\alpha\d s,
\end{equation}
where
\begin{eqnarray*}
 && a_\alpha(u,v)=\int_{\Omega^{2\pi}_H}\left[A_\zeta\grad u\cdot\grad\overline{v}-k^2 c_\zeta u\overline{v}\right]\d x-\int_{\Gamma^{2\pi}_H}\left(T^+_\alpha u\right)\overline{v}\d s;\\
 && b(u,v)=\int_{\Omega^{2\pi}_H}\left[(A_{\zeta_p}-A_\zeta)\grad u\cdot\grad\overline{v}-k^2(c_{\zeta_p}-c_\zeta)u\overline{v}\right]\d x;\\
&& F(\alpha,\cdot)=\frac{\partial(\J_\Omega u^i)(\alpha,\cdot)}{\partial x_2}-T_\alpha^+(\J_\Omega u^i)(\alpha,\cdot);\\
 && T^+_\alpha\phi=\i\sum_{j\in\Z}\sqrt{k^2-|j-\alpha|^2}\widehat{\phi}(j)e^{\i(j-\alpha)x_1},\,\phi=\sum_{j\in\Z}\widehat{\phi}(j)e^{\i(j-\alpha)x_1}.
\end{eqnarray*}

Following the arguments in \cite{Lechl2016,Lechl2017}, it is easy to prove that when the functions $\zeta$ and $\zeta_p$ are Lipschitz continuous, the variational problem \eqref{eq:var_bloch} is equivalent to \eqref{eq:sca_var}. When \eqref{eq:sca_var} has a unique solution of $\widetilde{H}_r^1(\Omega^p_H)$ for some $|r|<1$, the problem \eqref{eq:var_bloch} has a unique solution in $H_0^r(\W;\widetilde{H}^1_\alpha(D^{2\pi}))$. Moreover, {\color{lxl}if} the incident field $u^i\in H_r^2(\Omega^p_H)$ and the surfaces are $C^{2,1}$, then the solution belongs to the space $H_0^r(\W;\widetilde{H}^2_\alpha(D^{2\pi}))$. In \cite{Lechl2017}, a convergent numerical method based on \eqref{eq:var_bloch} has been proposed for the numerical solution, and a high order method has been proposed in \cite{Zhang2017e}.

\begin{remark}
 The information of the periodic function $\zeta$ is {\color{lxl}included} in $A_\zeta$ and $c_\zeta$, and the information of $p$ is {\color{lxl}included} in $A_{\zeta_p}$ and $c_{\zeta_p}$. During the iteration process, when $\zeta$ and $p$ are updated, the matrices $A_\zeta,\,A_p$ and the functions $c_\zeta,\,c_p$ are updated. Thus we do not need to change the meshes during {\color{lxl}this} process. 
\end{remark}

\section{Approximation of the scattering problems with periodic surfaces}

This section considers the difference between the scattered {\color{lxl}fields} with and without the local perturbation. Let $u_0$ be the total field with the same  incident field $u^i$ and periodic surface $\Gamma$, then $u_0$ satisfies the variational equation:
\begin{equation*}
  \int_{\Omega_H}\left[\grad u_0\cdot\grad\overline{v}-k^2 u_0\overline{v}\right]\d x-\int_{\Gamma_H}T^+\left(u_0\big|_{\Gamma_H}\right)\overline{v}\d s=\int_{\Gamma_H}f\overline{v}\d s.
\end{equation*}
From Theorem \ref{th:uni_sca}, {\color{lxl}if}  $f\in H_r^{-1/2}(\Gamma_H)$ for some $|r|<1$, then the solution $u_0\in \widetilde{H}_r^1(\Omega_H)$. In this paper, we assume that $r\in(0,1)$. 
As $u_0\in H_r^1(\Omega_H)$, there is a constant $C$ that does not depend on $u_0$ such that
\begin{equation*}
\left|u_0(x_1,x_2)\right|\leq C(1+x_1^2)^{-r/2-1/4}.
\end{equation*}

We apply the translation to the first variable, i.e., to replace $x_1$ by $x_1+2\pi L$ for some $L\in\Z\setminus\{0\}$, and let $u^i_L(x_1,x_2):=u^i(x_1+2\pi L,x_2)$ be the incident field. As the surface is $2\pi$-periodic, the total field with the incident field $u_L^i$, denoted by $u_0^L$, is actually the function $u_0(x_1+2\pi L,x_2)$. $u_0^L$ satisfies the following variational equation
\begin{equation*}
  \int_{\Omega_H}\left[\grad u_0^L\cdot\grad\overline{v}-k^2 u_0^L\overline{v}\right]\d x-\int_{\Gamma_H}T^+\left(u_0^L\big|_{\Gamma_H}\right)\overline{v}\d s=\int_{\Gamma_H}f_L\overline{v}\d s
\end{equation*}
with $f_L(x_1,x_2):=f(x_1+2\pi L,x_2)$ on $\Gamma_H$. As
\begin{equation*}
\left|u_0(x_1+2\pi L,x_2)\right|\leq C|2\pi L|^{-r-1/2},\quad  (x_1,x_2)\in\Omega^{2\pi}_H,
\end{equation*}
the following estimate holds:
\begin{equation*}
\left|u_0^L(x_1,x_2)\right|\leq C|2\pi L|^{-r-1/2},\quad  (x_1,x_2)\in\Omega^{2\pi}_H.
\end{equation*}

Let $u^L$ be the solution of \eqref{eq:sca_var} with $f$ be replaced by $f_L$. Similar to the previous section, we {\color{lxl}can define} a diffeomorphism $\Phi_p$ that maps $\Omega^p_H$ to $\Omega_H$ and $\Phi_p-I_2$ is supported in $\Omega_H\cap W\times\R$. Let $u_T^L:=u^L\circ \Phi_p$, {\color{lxl}it is} easily checked that $u^L_T$ satisfies 
\begin{equation}\label{eq:var_T}
  \int_{\Omega_H}\left[A_p\grad u_T^L\cdot\grad\overline{v}-k^2 c_p u_T^L\overline{v}\right]\d x-\int_{\Gamma_H}T^+\left(u_T^L\big|_{\Gamma_H}\right)\overline{v}\d s=\int_{\Gamma_H}f_L\overline{v}\d s,
\end{equation}
where 
\begin{eqnarray*}
 && A_p(x)=\left|\det\grad\Phi_p(x)\right|\left[\left(\grad\Phi_p(x)\right)^{-1}\left(\left(\grad\Phi_p(x)\right)^{-1}\right)^\top\right]\in L^\infty(D,\R^{2\times 2});\\
 &&c_p(x)=\left|\det\grad\Phi_p(x)\right|\in L^\infty(D).
\end{eqnarray*}
Moreover, ${\rm supp}(A_p-I_2),\,{\rm supp}(c_p-1)\subset\Omega_H\cap W\times\R(:=\Omega^{2\pi}_H)$. Then the difference $u_d^L:=u_T^L-u_0^L$ satisfies the following variational equation, i.e.,
\begin{equation}\label{eq:var_diff}
\int_{\Omega_H}\left[A_p\grad u_d^L\cdot\grad\overline{v}-k^2 c_p u_d^L\overline{v}\right]-\int_{\Gamma_H}T^+\left(u^L_d\big|_{\Gamma_H}\right)\overline{v}\d s=\widetilde{b}(u_0^L,v)
\end{equation}
for any $v\in H^1(\Omega_H)$ with compact support, where
\begin{equation}
 \widetilde{b}(u_0^L,v)=\int_{\Omega^{2\pi}_H}\left[(I_2-A_p)\grad u_0^L\cdot\grad\overline{v}-k^2 (1-c_p) u_0^L\overline{v}\right]\d x.
\end{equation}
From the representation of $\widetilde{b}(\cdot,\cdot)$, it is a bounded sesquilinear form satisfies
\begin{equation*}
\left|\widetilde{b}(u,v)\right|\leq C\|u\|_{H^1(\Omega^{2\pi}_H)}\|v\|_{H^1(\Omega^{2\pi}_H)},
\end{equation*}
where $C$ is the constant depends only on $\zeta$ and $\zeta_p$. 
Thus the right hand side of \eqref{eq:var_diff} satisfies
\begin{equation*}
\left|\widetilde{b}(u_0^L,v)\right|\leq C\left\|u_0^L\right\|_{H^1(\Omega^{2\pi}_H)}\|v\|_{H^1(\Omega^{2\pi}_H)}\leq C|2\pi L|^{-r-1/2}\|v\|_{H^1(\Omega_H)}.
\end{equation*}
From the equivalence between \eqref{eq:sca_var} and \eqref{eq:var_trans}, the equation \eqref{eq:var_diff} is uniquely solvable in $H^1(\Omega_H)$ when the right hand side is a antilinear functional on $H^1(\Omega_H)$. Thus
\begin{equation*}
\left\|u_T^L-u_0^L\right\|_{H^1(\Omega_H)}\leq C|2\pi L|^{-r-1/2}.
\end{equation*}

{\color{lxl}Based on the above analysis,  the total field  $u_0^L$ is a good approximation of $u^L_T$ if $L$ is sufficiently large}. Especially, let $\sigma$ be the noise level of the measured data. {\color{lxl}When} $L\in\Z$ has a large enough absolute value such that $C|2\pi L|^{-r-1/2}<\delta$, $u^L_T$ could be treated as the ``exact solution'' of the non-perturbed periodic surface with the incident field $u_T^L$. Let
\begin{equation*}
 \widetilde{u}^L_T(x_1,x_2):=u^L_T(x_1-2\pi L,x_2),
\end{equation*}
then $\widetilde{u}^L_T$ is a good approximation of $u_0$. In this case, the solution $\widetilde{u}^L_T$ could be applied in the inverse problems to reconstruct the periodic surface.

\section{Inverse Problem and the Newton-CG Method}

The aim of the inverse problem is to reconstruct the unknown function $\zeta_p$ from the measured scattered data. The measured scattered field $U$ on $\Gamma_H$ is defined as
\begin{equation}
 U:=u^s\big|_{\Gamma_H}+\sigma
\end{equation}
where $\sigma$ is some noise added to the scattered data.

In the following, we always assume that  $\zeta\in C^{2,1}(\R)$ is a $2\pi$-periodic function and $p\in C^{2,1}(\R)$ is a function that is compactly supported in $W+2\pi J$ for some $J\in\Z$.
\begin{remark}
 For the inverse problem, $J$ is an unknown integer and one task for the inverse problem is to find out the exact value of $J$. As is explained later, the integer $J$ could be found out by a sampling method (see \cite{Liu2018}). Thus in this section, we treat it as a known one. For any $J\neq 0$, we can simply apply the translation $x\mapsto x-2\pi J$ to move the perturbation to the center of the domain (i.e., $J=0$). Thus for simplicity, we still assume that $J=0$ {\color{lxl} in this section}.
\end{remark}

Define the spaces
\begin{eqnarray*}
&& X:=\{\zeta\in C^{2,1}(\R):\,\zeta\text{ is }2\pi-\text{periodic}\};\\
&&  Y:=\{p\in C^{2,1}(\R):\,{\rm supp}(p)\subset W\}.
\end{eqnarray*}
In the following, we assume that $(\zeta,p)\in X\times Y$ and $\zeta_p:=\zeta+p$. {\color{lxl}The} inverse problem is to find out {\color{lxl}$(\zeta,p)\in X\times Y$} such that the scattered field {\color{lxl}corresponding} to $(\zeta,p)$ is the best approximation of $U$. 

\subsection{Scattering operator and its properties with respect to rough surfaces}

We recall the inverse scattering problems from rough surfaces introduced in \cite{Chand2002a}. Let $BC^{1,1}(\R)$ be the space of bounded, Lipschitz continuous function. 

\begin{remark}
 It is easily checked that $X,\,Y\subset BC^{1,1}(\R)$.
\end{remark}

Suppose $f\in BC^{1,1}(\R)$ and the surface $\Gamma_f$ is defined by $f$. We can also define the domain $\Omega_f$ by the domain above $\Gamma_f$, and $\Omega^f_H$ by the domain between $\Gamma_f$ and $\Gamma_H$, where $H$ is a real number that is larger than $\|f\|_\infty$. Given an incident field $u^i$, we define the following scattering operator
\begin{eqnarray*}
  S:\quad BC^{1,1}(\R)&\rightarrow & L^2(\Gamma_H)\\
 f&\mapsto & u^s\big|_{\Gamma_H}.
\end{eqnarray*}

Then the inverse problem {\color{lxl}can be} written as the optimization problem, i.e., to find $f\in BC^{1,1}(\R)$ such that
\begin{equation}
 \|S(f)-U\|^2_{L^2(\Gamma_H)}=\min_{f^*\in BC^{1,1}(\R)}\|S(f^*)-U\|^2_{L^2(\Gamma_H)}.
\end{equation}
Let
\begin{equation}
 F(f):= \|S(f)-U\|^2_{L^2(\Gamma_H)},
\end{equation}
then the inverse problem is to find out the minimizer of the functional $F$ in the domain $BC^{1,1}(\R)$. To solve the minimization problem, we have to study the properties of the scattering operator $S$ first.

\begin{theorem}
 The operator $S$ is differentialble, and its derivative $DS$ is represented as
 \begin{align}
   DS:\, BC^{1,1}(\R) &\rightarrow L^2(\Gamma_H)\\ 
    h &\mapsto u'\big|_{\Gamma_H},
      \end{align}
 where $u'\in H_r^1(\Omega^f_H)$ satisfies
 \begin{eqnarray}\label{eq:sca_diff1}
  && \Delta u'+k^2 u'=0\quad\text{ in }\Omega^fH;\\\label{eq:sca_diff2}
  && u'=-\frac{\partial u}{\partial x_2}h\quad\text{ on }\Gamma_f;\\\label{eq:sca_diff3}
  && \frac{\partial u'}{\partial x_2}=T^+u'\quad\text{ on }\Gamma_H.
 \end{eqnarray}
Here $u$ is the total field of the scattering problem \eqref{eq:sca1}-\eqref{eq:sca3}.
\end{theorem}
For the proof of this theorem we refer to \cite{Kirsc1993b,Chand2002a}.

For the Newton's method, we also need the adjoint operator of the Fr{\'e}chet derivative $DS$, which is explained in the following Theorem.

\begin{theorem}
 The adjoint operator of $DS(f)$, denoted by $\left[DS(f)\right]^*$ is given by
 \begin{equation}
  \left[DS(f)\right]^*\phi=-\Re\left[\frac{\partial\overline{u}}{\partial\nu}\frac{\partial \overline{z}}{\partial\nu}\right]\nu_2,
 \end{equation}
where $\nu$ is the normal derivative upwards, $u$ is the total field and $z$ satisfies
\begin{eqnarray}\label{eq:sca_adj1}
 && \Delta z+k^2 z=0\quad\text{ in }\Omega_H^f;\\\label{eq:sca_adj2}
 && z=0\quad\text{ on }\Gamma_f;\\\label{eq:sca_adj3}
 && \frac{\partial z}{\partial x_2}-T^+z=\overline{\phi}\quad\text{ on }\Gamma_H.
\end{eqnarray}
\end{theorem}

\begin{remark}
During the iteration steps, the problems \eqref{eq:sca_diff1}-\eqref{eq:sca_diff3} and \eqref{eq:sca_adj1}-\eqref{eq:sca_adj3} will be solved several times. We can always apply the method introduced in Section 2.2 to transform the problems first into the one defined in the unbounded rectangle $D$ by the transform $\Phi_{\zeta_p}$, and then apply the Floquet-Bloch transform to obtain the new problem defined in the bounded domain $\W\times D^{2\pi}$. For details of the solution of \eqref{eq:sca_adj1}-\eqref{eq:sca_adj3} we refer to Remark 12 in \cite{Zhang2017c}. 
\end{remark}

\subsection{Discretization for locally perturbed periodic surfaces}

Let $\{\phi_1,\phi_2,\dots,\phi_M,\dots\}$ be a basis in the space $X$ and $\{\psi_1,\psi_2,\dots,\psi_N,\dots\}$ be a basis in the space $Y$. For the positive integers $M$ and $N$,  define the finite-dimensional subspaces of $X$ and $Y$ by:
\begin{equation*}
 X_M:={\rm span}\{\phi_1,\dots,\phi_M\}\subset X\quad{ and }\quad Y_N:={\rm span}\{\psi_1,\dots,\psi_N\}\subset Y.
\end{equation*}
For the coefficients ${\bm C^M}=(c_1^M,\dots,c_M^M)\in\R^M$ and ${\bm D^N}=(d_1^N,\dots,d_N^N)\in\R^N$, then the elements $\zeta_M\in X_M$ and $p_N\in Y_N$ could be written as
\begin{equation*}
 \zeta_M(t)=\sum_{m=1}^M c_m^M\phi_m(t),\quad p_N(t)=\sum_{n=1}^N d_n^N\psi_n(t).
\end{equation*}
For the function $\zeta\in X$, there is a ${\bm C^M}\in\R^M$ such that $\zeta_M$ is the approximation in $X_M$ of $\zeta$. The argument also holds for $p_N\in Y_N$ and $p\in Y$.

Define the operators $A$ and $B$ by
\begin{equation*}
 \begin{aligned}
  A:\,\R^M &\rightarrow X_M\\
  {\bm C^M} &\mapsto \zeta_M
 \end{aligned};\quad
  \begin{aligned}
  B:\,\R^N &\rightarrow Y_N\\
  {\bm D^N} &\mapsto p_N.
 \end{aligned}
\end{equation*}
Then define the operator
\begin{equation*}
 \begin{aligned}
P:\, \R^M\times\R^N & \rightarrow L^2(\Gamma_H)\\
\left({\bm C^M},\,{\bm D^N}\right)&\mapsto S\circ\left(A({\bm C^M})+B({\bm D^N})\right), 
\end{aligned}
\end{equation*}
which maps the coefficients of both the periodic function and the local perturbation to the scattered field. Then we can define the functional $F$ in the finite dimensional space $\R^M\times\R^N$ by
\begin{equation}
 F({\bm C^M},{\bm D^N}):=\|P({\bm C^M},{\bm D^N})-U\|^2_{L^2(\Gamma_H)},
\end{equation}
and the inverse problem is formulated by the following finite dimensional problem:\\
\noindent
{\bf Discrete Inverse Problem: } to find $({\bm C^M},{\bm D^N})\in\R^M\times\R^N$ such that 
\begin{equation}\label{eq:dis_inv}
 F({\bm C^M},{\bm D^N})=\min_{({\bm C^M_*},{\bm D^N_*})\in\R^M\times\R^N}F({\bm C^M_*},{\bm D^N_*}).
\end{equation}

We apply the Newton-CG method to solve the descritized inverse problem. The linearized equation is
\begin{equation}
 P({\bm C^M},{\bm D^N})+(DP)({\bm C^M},{\bm D^N})(\delta{\bm C^M},\delta{\bm D^N})=U,
\end{equation}
where $\delta{\bm C^M}=(\delta c_1^M,\dots,\delta c_M^M)\in\R^M$ and $\delta{\bm D^N}=(\delta d_1^N,\dots,\delta d_N^N)\in\R^N$, $(DP)({\bm C^M},{\bm D^N})$ is the Fr{\'e}chet derivative of $P$ at $({\bm C^M},{\bm D^N})$.  Define 
\begin{eqnarray*}
&&M_A({C^M})(\delta{\bm C^M}):=(DP)({\bm C^M},{\bm D^N})(\delta{\bm C^M},{\bm 0});\\
 &&M_B({D^N})(\delta{\bm D^N}):=(DP)({\bm C^M},{\bm D^N})({\bm 0},\delta{\bm D^N}),
\end{eqnarray*} 
then the linearized equation is written as
\begin{equation}
 P({\bm C^M},{\bm D^N})+M_A({\bm C^M})(\delta{\bm C^M})+M_B({\bm D^N})(\delta{\bm D^N})=U.
\end{equation}

First, we have to calculate the derivative of $P$. As an operator defined in the finite dimensional space $\R^M\times\R^N$, from direct calculation,
\begin{equation*}
\frac{\partial P}{\partial c_m^M}=(DS)(A({\bm C^M})+B({\bm D^N}))\phi_m;\quad
\frac{\partial P}{\partial d_n^N}=(DS)(A({\bm C^M})+B({\bm D^N}))\psi_n.
\end{equation*}
Thus
\begin{eqnarray*}
 && M_A({\bm C^M})(\delta{\bm C^M})=\sum_{m=1}^M \delta c_m^M \frac{\partial P}{\partial c_m^M}= (DS)(A({\bm C^M})+B({\bm D^N}))\left[\sum_{m=1}^M \delta c_m^M\phi_m\right];\\
 && M_B({\bm D^N})(\delta{\bm D^N})=\sum_{n=1}^N \delta d_n^N \frac{\partial P}{\partial d_n^N}= (DS)(A({\bm C^M})+B({\bm D^N}))\left[\sum_{n=1}^N \delta d_n^N\psi_n\right].
\end{eqnarray*}

Given any $\delta {\bm C^M}\in\R^M$ and $\phi\in L^2(\Gamma_H)$,
\begin{equation*}
 \begin{aligned}
  \left(\delta{\bm C^M},M_A^*({\bm C^M })\phi\right)&=\left((M_A)({\bm C^M})(\delta{\bm C^M}),\phi\right)\\
  &=\left((DS)(A({\bm C^M})+B({\bm D^N}))\left[\sum_{m=1}^M \delta c_m^M\phi_m\right],\phi\right)\\
  &=\sum_{m=1}^M \delta c_m^M\left(\phi_m,\left[(DS)(A({\bm C^M})+B({\bm D^N}))\right]^*\phi\right).
 \end{aligned}
\end{equation*}
Let $Q=\left[(DS)(A({\bm C^M})+B({\bm D^N}))\right]^*$, then
\begin{equation}
 M_A^*({\bm C^M })\phi=\left((\phi_1,Q\phi),\dots,(\phi_M,Q\phi)\right).
\end{equation}
Similarly, we can also get
\begin{equation}
 M_B^*({\bm D^N })\phi=\left((\psi_1,Q\phi),\dots,(\psi_N,Q\phi)\right).
\end{equation}

In the numerical implementation, we solve the discrete inverse problem separately, i.e., first fix ${\bm D^N}$ and solve the minimization problem \eqref{eq:dis_inv} {\color{lxl}to} find out the solution ${\bm C^M}$.{\color{lxl}Then we} fix ${\bm C^M}$ and solve the problem with respect to ${\bm D^N}$.

To solve the minimization problems we apply the Newton-CG method. To minimize the function $F({\bm C^M},{\bm D^N})$ with fixed ${\bm D^N}$, we apply the following Newton-CG method.
\begin{algorithm}[H]\caption{Newton-CG Method -- Part I}
Input: Data $U$; $\epsilon>0$; $j=0$; fixed ${\bm D^N}\in\R^N$.\\
Initialization: ${\bm C^M_0}\in\R^M$.
\begin{algorithmic}[1] 
\While{$\|P({\bm C^M_j},{\bm D^N})-U\|_{L^2(\Gamma_H)}>\epsilon\|U\|_{L^2(\Gamma_H)}$}
\\\quad $CGNE$ iteration scheme to solve $M_A({\bm C^M_j})({\bm H^M})=U-P({\bm C^M_j},{\bm D^N})$
\\\quad
${\bm C^M_{j+1}}={\bm C^M_j}+{\bm H^M}$;\\
\quad $j=j+1$;
\EndWhile
\end{algorithmic}
\label{alg1}
\end{algorithm}
Similarly, we can also minimize the function $F({\bm C^M},{\bm D^N})$ with fixed ${\bm C^M}$ by the following algorithm:
\begin{algorithm}[H]\caption{Newton-CG Method -- Part II}
Input: Data $U$; $\epsilon>0$; $j=0$; fixed ${\bm C^M}\in\R^M$.\\
Initialization: ${\bm D^N_0}\in\R^N$.
\begin{algorithmic}[1] 
\While{$\|P({\bm C^M},{\bm D^N_j})-U\|_{L^2(\Gamma_H)}>\epsilon\|U\|_{L^2(\Gamma_H)}$}
\\\quad $CGNE$ iteration scheme to solve $M_B({\bm D^N_j})({\bm H^N})=U-P({\bm C^M},{\bm D^N_j})$
\\\quad
${\bm D^N_{j+1}}={\bm D^N_j}+{\bm H^N}$;\\
\quad $j=j+1$;
\EndWhile
\end{algorithmic}
\label{alg2}
\end{algorithm}

\section{Numerical implementation}
\subsection{Sampling method}

{\color{xl}
In this section, we use the sampling method introduced in \cite{Liu2018} to give an initial guess of 
the perturbed periodic surface, especially for the first term ${\bm c_1^0}$ of ${\bm C^0}$ and the integer $J$ of the perturbation.

Suppose that location $y$ of incident point sources $u^i(x,y)$ is on a horizontal line $\Gamma_H:=\{(x_1, H)~|~x_1\in\R\}$ above the surface, we measure scattered Cauchy data $(u^s, \partial_{\nu} u^s)$ generated by these point sources and the perturbed periodic surface on $\Gamma_H$. Here, $\partial_{\nu} u^s$ denotes the normal derivative of $u^s$ on $\Gamma_H$ with the direction $(0,1)$.

We introduce the following imaging function
\begin{align}\nonumber
&I(z)=\int_{\Gamma_H}\left|\int_{\Gamma_H}\left({\partial_{\nu(x)} u^s(x,y)}\ol{\varPhi_{k}(x,z)}
-u^s(x,y){\partial_{\nu(x)} \ol{\varPhi_{k}(x,z)}}\right)ds(x)\right.\\ \label{eq18}
&\qquad\qquad\left.-\frac{i}{4\pi}\int_{\mathbb{S}_-}e^{ik\hat{x}\cdot(y'-z')}ds(\hat{x})
\right|^2ds(y),
\end{align}
where $y'=(y_1,-y_2)$ and $z'=(z_1,-z_2)$. From the analysis in \cite{Liu2018}, we can expect that the imaging function $I(z)$ takes a large value when $z\in\Gamma_p$ and decays
as $z$ moves away from $\Gamma_p$. In this way, we give an initial guess of the perturbed surface.

In numerical computation, we choose $2P+1$ incident point sources which are located at $y_j=(jh_{inc},H),$ $j=-P,...,0,...P,$ here $h_{inc}$ is a fixed interval between two adjacent points.
The measurement line $\Gamma_H$ is truncated to be ${\Gamma_{H,A}}:=\{x\in\Gamma_H~|~|x_1|<A\}$ which will be discretized uniformly into $2Q$ subintervals so the step size is $h_{mea}=A/Q$. In addition, the lower-half circle $\mathbb{S}_-$ in the second integral
in (\ref{eq18}) will also be uniformly discretized into $R$ grids with the step size $\Delta\theta=\pi/R$. Then for each sampling point $z$ we get the following discrete form of (\ref{eq18})
\begin{align}\nonumber
&I_A(z)=\sum_{j=-P}^{P}\left|h_{mea}\sum_{i=0}^{2Q}
\left({\partial_{\nu(x)} u^s(x_i,y_j)}\ol{\varPhi_{k}(x_i,z)}
-u^s(x_i,y_j){\partial_{\nu(x)}\ol{\varPhi_{k}(x_i,z)}}\right)\right.\\ \label{eq188}
&\qquad\qquad\left.-\frac{i\Delta\theta}{4\pi}
\sum_{k=0}^{R}e^{ikd_k\cdot(y_j'-z')}\right|^2.
\end{align}
Here, the measurement points are denoted by $x_i=(-A+ih,H),$ $i=0,1,...,2Q,$
and the normal directions are denoted by $d_k=(\sin(-\pi+k\Delta\theta),\cos(-\pi+k\Delta\theta)),$ $k=0,1,...,R.$

Suppose the sampling area is a rectangle denoted by $[a, b]\times[c, d]$. We set the numbers of sampling points in $x_1$-direction and $x_2$-direction to be $M_1$ and $M_2$, respectively. Then by (\ref{eq188}), we get the indicator matrix $\left\{I_A(z_{ij})\right\}_{M_1\times M_2}$. For each $j$th-row of this matrix, we figure out the element with the largest value $I_A$ and denote the corresponding index by ${{max}}_{j}$.  The initial guess for the first term ${\bm c^0_1}$ of $C^1$ can be deduced by the following formula
\begin{equation}\label{initial}
{\bm c^0_1}= c+(d-c)\frac{1}{M_1 M_2}\sum_{j=1}^{j=M_1} max_j.
\end{equation}
}

\subsection{Iteration method}

From the last subsection, we have already decided the integer $J$. By translation on the first variable, i.e., to let $x_1$ be replaced with $x_1-2\pi J$, {\color{lxl}the perturbation} is moved to  $W$ (i.e., $J=0$), then the method in Section 4 could be applied to the reconstruction. 

From Section 3, for an incident field $u^i\in H_r^1(\Omega^p_H)$ for some $r\in(0,1)$, the measured  {\color{lxl}data $U_L$} with the incident field $u^i(\cdot+2\pi L,\cdot)$ for $L\in\Z\setminus\{0\}$,  {\color{lxl}could be} applied to reconstruct the periodic function $\zeta$. The measured data with incident field $u^i$, denoted by $U_0$, is then applied to reconstruct the local perturbation $p$. So we conclude the algorithm for the inverse scattering problem.

\begin{algorithm}[H]\label{alg}
\caption{Numerical Method for the Inverse Problem}
Input: Cauchy date $(u^s_j,\partial_\nu u^s_j)$  generated by point sources located at $x_j$;\\
Given: Domain $D$, $\mathcal{M}$ is a regular mesh for $D$.
\begin{enumerate}
\item Decide $J$ and $c_0^0$ from the sampling method. Move the perturbation to the center by $x_1\mapsto x_1-2\pi J$.\\
Generate the measured date $U_0$ and $U_L$ with incident fields $u^i(x_1,x_2)$ and $u^i(x_1+2\pi L,x_2)$.\\
Set the initial guess:  ${\bm C_0}:=({\bm c_0^0},0,\dots,0)$, ${\bm D_0}={\bm 0}$.
\item Solve the minimization problem for the fixed ${\bm D_0}$ by Algorithm \ref{alg1}:
\begin{equation*}
F({\bm C},{\bm D_0})=\|P({\bm C},{\bm D_0})-U_L\|\rightarrow\min.
\end{equation*}
\item Solving the minimization problem for the fixed ${\bm C}$ by Algorithm \ref{alg2}:
\begin{equation*}
F({\bm C},{\bm D})=\|P({\bm C},{\bm D})-U_0\|\rightarrow\min.
\end{equation*}
\end{enumerate} 
Then  $({\bm C},{\bm D})$ is the final result of the numerical scheme.
\end{algorithm}

\section{Numerical results}

\begin{figure}[H]
\centering
\begin{tabular}{c c}
\includegraphics[width=0.45\textwidth]{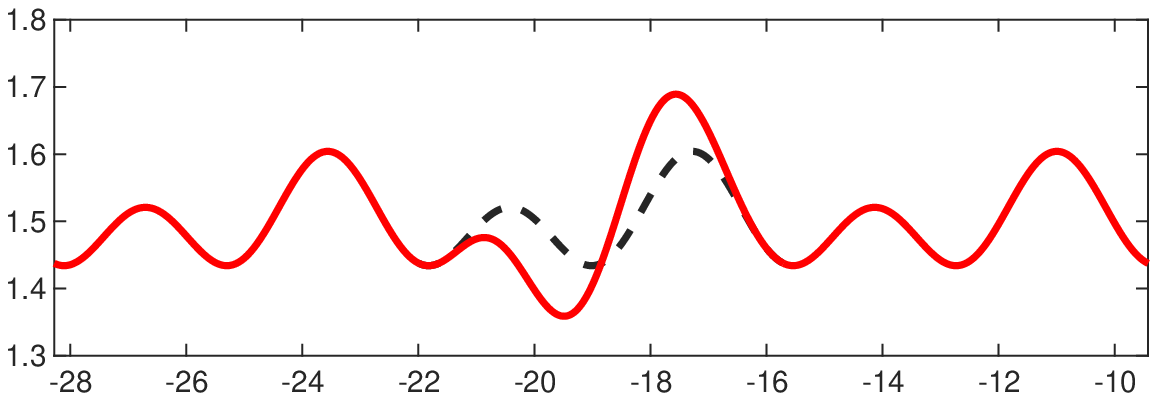} 
& \includegraphics[width=0.45\textwidth]{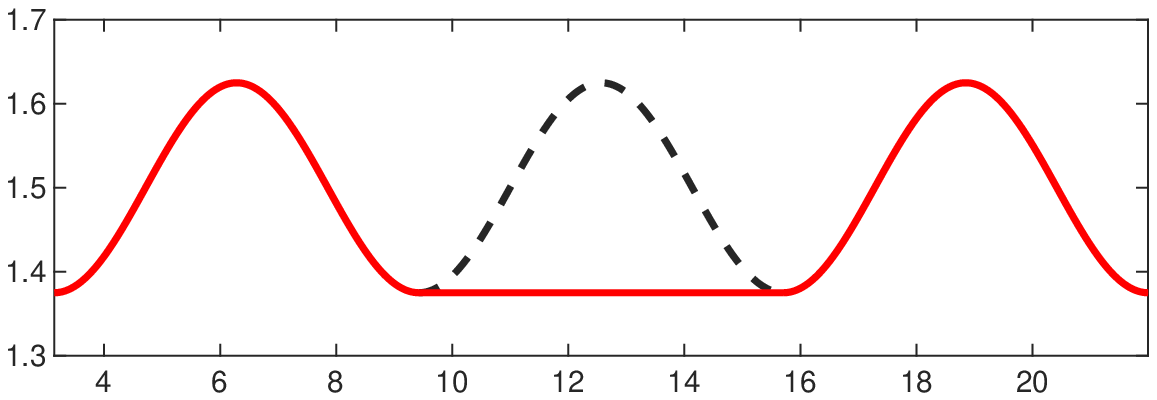}\\[-0cm]
(a) & (b) 
\end{tabular}%
\caption{(a): the first surface; (b): the second surface. }
\label{fig:surface}
\end{figure}

In this section, we present two examples for the numerical method. We define two different periodic surfaces and local perturbations:
\begin{eqnarray*}
 && \zeta_1(t)=1.5+\frac{\sin t}{24}-\frac{\cos 2t}{16};\\
 && \zeta_2(t)=1.5+\frac{\cos t}{8};\\
 && p_1(t)=0.00025 ((t+6\pi)^2-9)^3\sin\left(\frac{\pi(t+3)}{3}\right)\mathcal{X}_{[-3-6\pi,3-6\pi]}(t);\\
 && p_2(t)=-\frac{1+\cos t}{8}\mathcal{X}_{[-3+4\pi,3+4\pi]}(t).
\end{eqnarray*}
We apply Algorithm \ref{alg} to the following two examples (see Figure \ref{fig:surface}):\\

\vspace{0.15cm}
\noindent
{\bf Example 1.} The periodic surface $\Gamma$ is defined by $\zeta_1$ and the local perturbation is defined by $p_1$;\\
\vspace{0.1cm}

\noindent
{\bf Example 2.} The periodic surface $\Gamma$ is defined by $\zeta_2$ and the local perturbation is defined by $p_2$.
\vspace{0.15cm}

For both the incident point sources and Herglotz wave functions, the scattered data are collected on $\Gamma_{A,H}$ with $A=25\pi, H=3$
and it is divided into $2Q=1500$ subintervals with the step length $h_{mea}=\pi/300$. Let $u_{s}$ be the scattered data (either the scattered field or its normal derivative) on $\Gamma_{A,H}$,  {\color{lxl}then} the measured data is defined as:
\begin{equation*}
U_{meas}:=u_s+\sigma\max(u_s){\rm randn},
\end{equation*}
where $\sigma=5\%$ is the noise level and randn presents random numbers from the standard normal distribution.

\subsection{Sampling method}
{\color{xl}

For the sampling method, we choose the sampling area to be a rectangle as $[-20\pi,20\pi]\times [1.2,1.9]$. The number of sampling points in $x_1$-direction and $x_2$-direction are set to be $M_1=1600$ and $M_2=400$, respectively. For the first surface, we put 41 incident point sources at $y_j=(j\pi, 3)$ with $j=-20,-19,\dots,20.$  For the second surface, we put 21 incident point sources at $y_j=(2j\pi, 3)$ with $j=-10,-9,\dots,10$. The wavenumber is chosen to be $k=3$ for both examples. 

Use the indicator function introduced in (\ref{eq188}), we can get a rough reconstruction of the original perturbed periodic surfaces in Figure \ref{fig:recon1} and \ref{fig:recon2}. Note that the red dash lines in (c) are boundaries of the periodic cells.
In each figure, we first present the profile of the original surface. Then the reconstructed result is given directly by the indicator function $I_{A}(z)$. Finally, in order to give the initail guess of ${\bm c^0_1}$ and the integer $J$ of the perturbation, we try to find out the points $z_{max}$ which get the largest value $I_A(z)$ in each vertical line and plot them in the last position of each figure.

\begin{figure}[H]
\centering
\includegraphics[width=\textwidth]{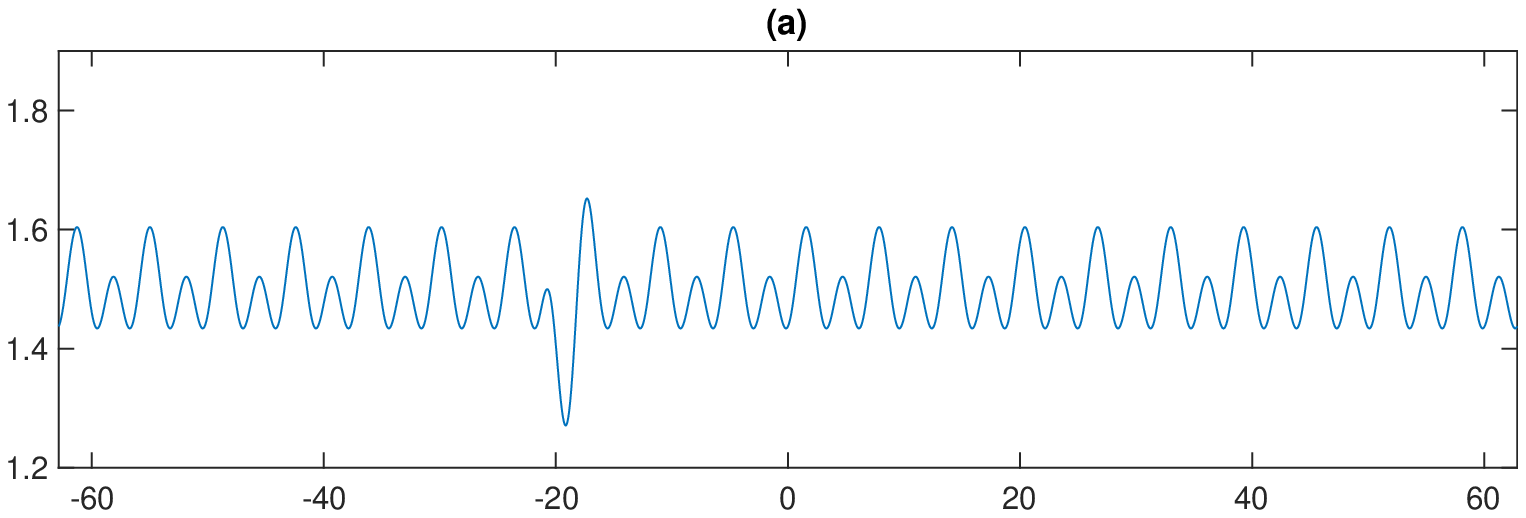} 
\includegraphics[width=\textwidth]{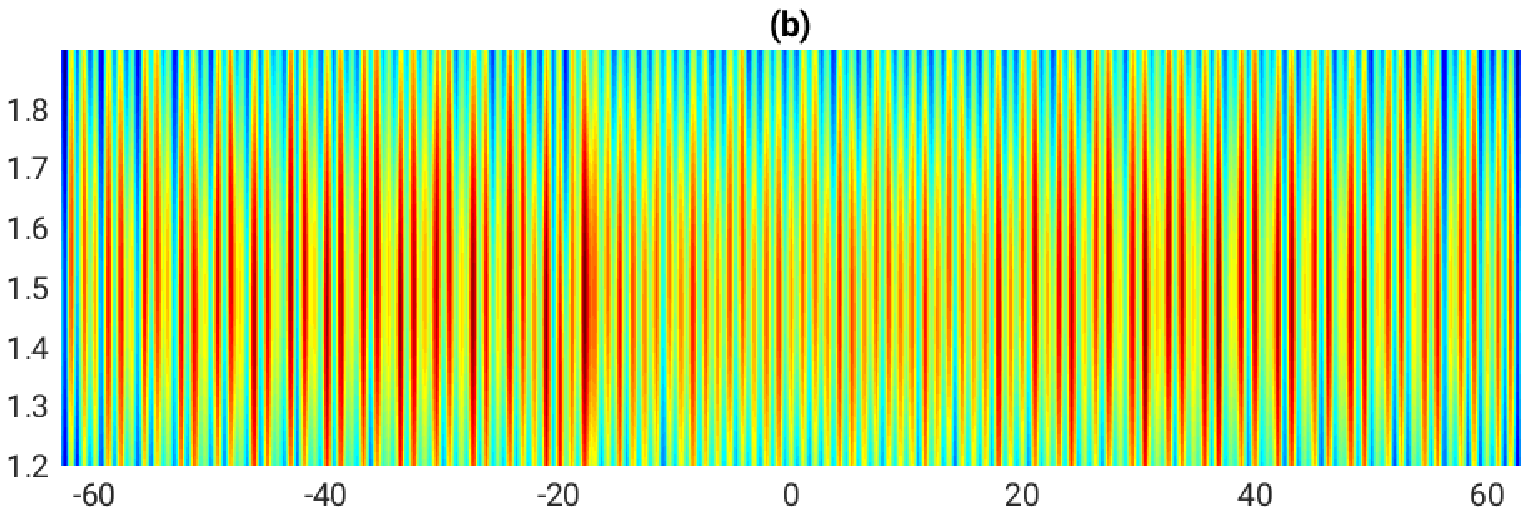} 
\includegraphics[width=\textwidth]{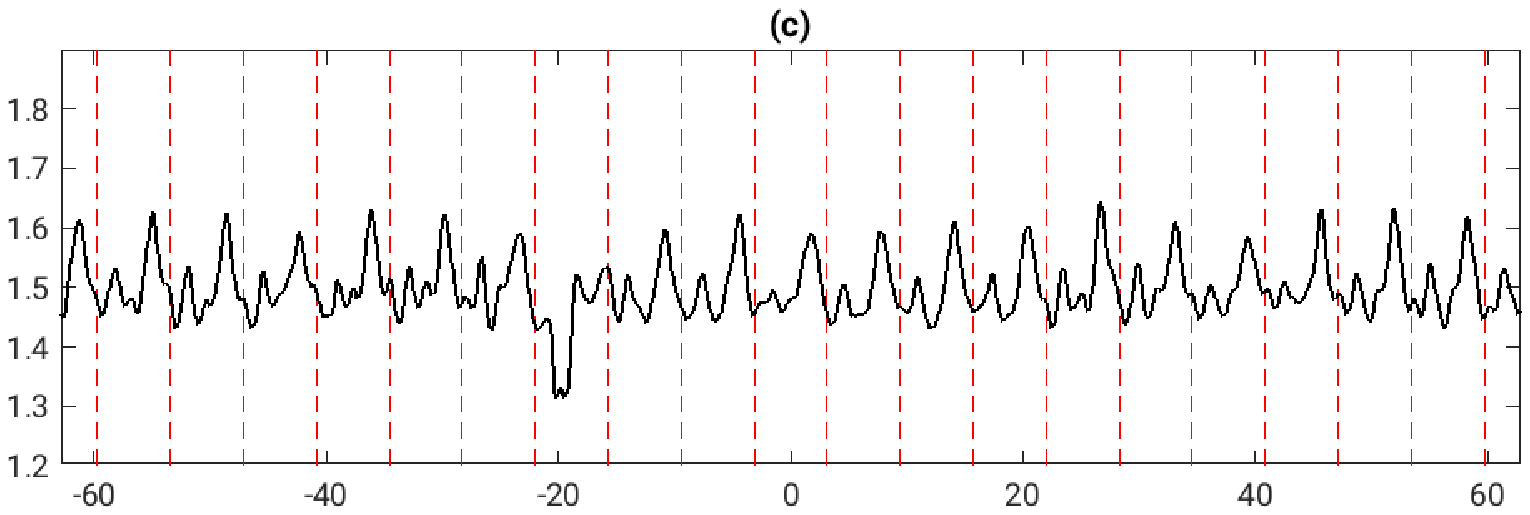} 
\caption{(a): the first surface; (b)\&(c): the reconstructions. }
\label{fig:recon1}
\end{figure}

\begin{figure}[H]
\centering
\includegraphics[width=\textwidth]{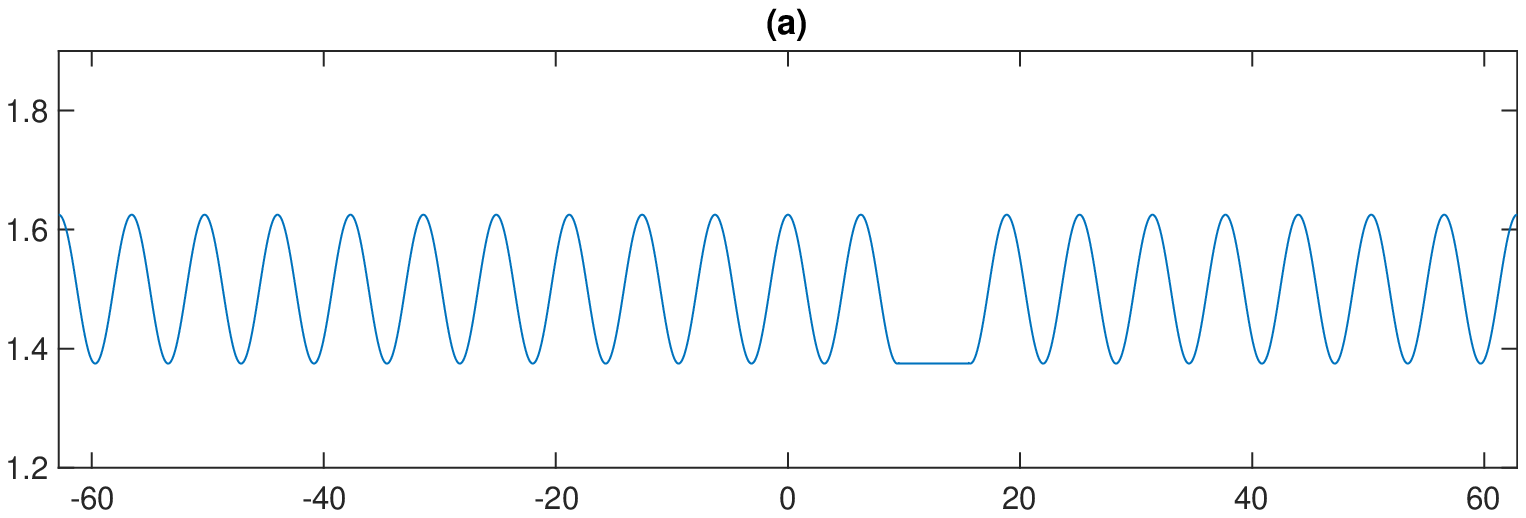} 
\includegraphics[width=\textwidth]{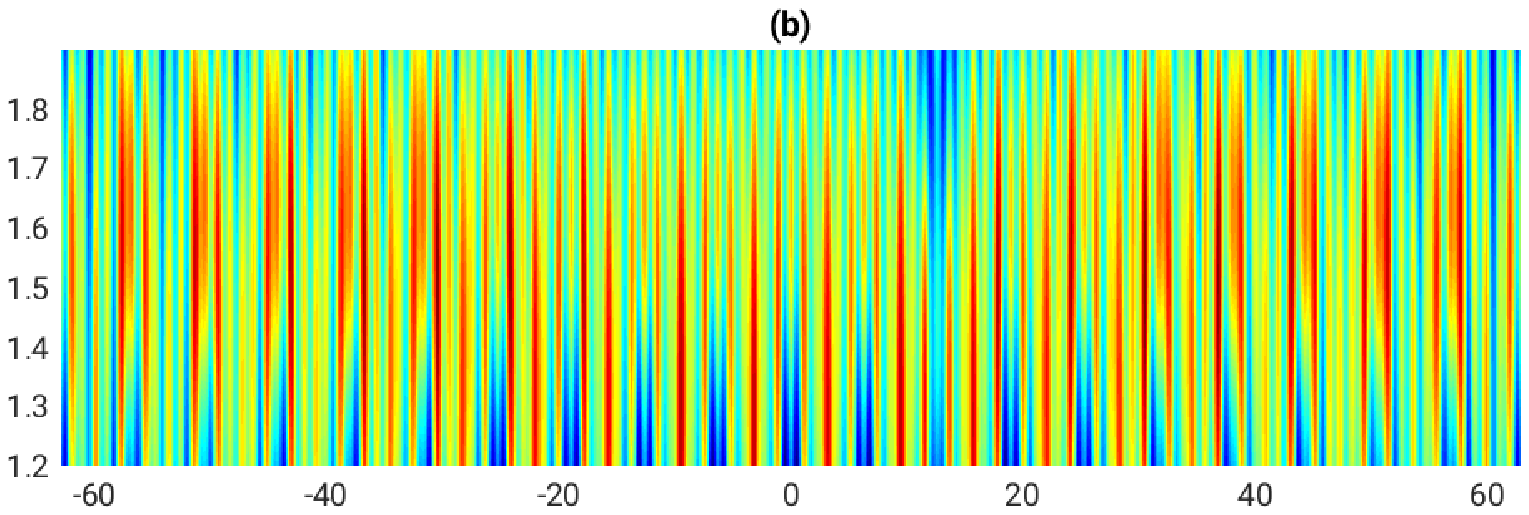} 
\includegraphics[width=\textwidth]{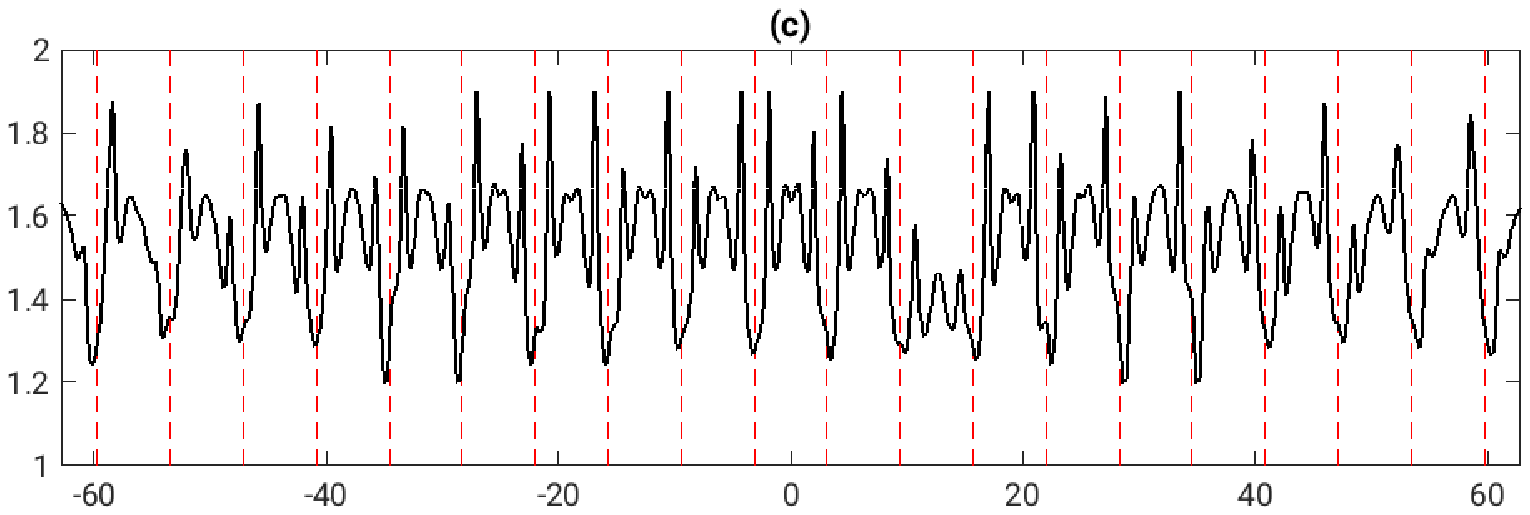} 
\caption{(a): the first surface; (b)\&(c): the reconstructions. }
\label{fig:recon2}
\end{figure}

By the end of the sampling step, we give out the value $J$ and ${\bm c^0_1}$. 
Roughly speaking, $J$ represents the location of the perturbation while ${\bm c^0_1}$ gives the vertical location of the periodic surface. 
From Figure \ref{fig:recon1} and \ref{fig:recon2}, the locations of the perturbations are easily obtained, i.e., $J=-3$ for Example 1 and $J=2$ for Example 2.
The initial guess of ${\bm c^0_1}$ is computed due to (\ref{initial}). By staightward calcultaions, we get ${\bm c^0_1}=1.4987$ and ${\bm c^0_1}=1.5216$. Both of these two results are very good approximations of the constant terms of both $\zeta_1$ and $\zeta_2$.
}

\subsection{Newton's method}

\begin{figure}[H]
\centering
 \includegraphics[width=16cm]{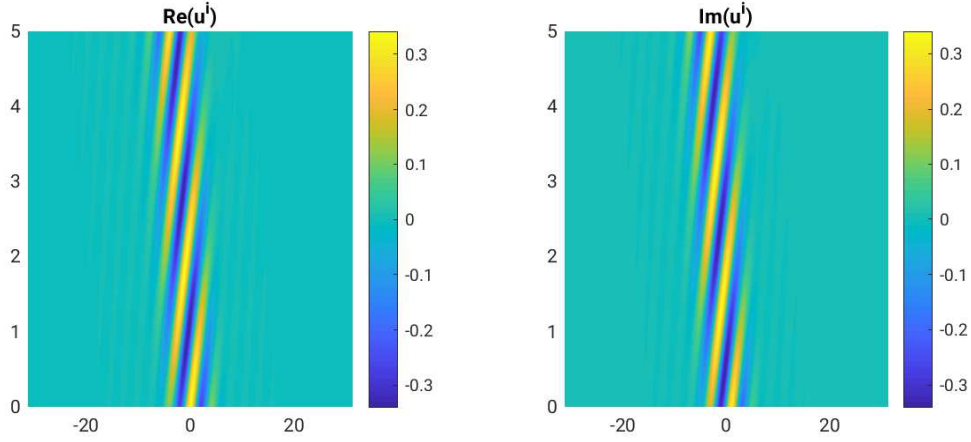}
 \caption{Real- and imaginary-part of the incident field $u^i$.}
 \label{fig:incident}
\end{figure}

For the Newton's method, {\color{lxl}the Herglotz wave function is} applied as the incident field (see Figure \ref{fig:incident}), i.e.,
\begin{equation*}
 u^i(x_1,x_2)=\int_{-\pi/2}^{\pi/2}\exp\left(\i k (x_1\sin t-x_2 \cos t)\right)g(t)\d t,
\end{equation*}
where 
\begin{equation*}
 g(t)=2^{12} t^6(1-t)^6\mathcal{X}_{[0,1]}(t).
\end{equation*}

\begin{remark}
We could not use the point source as the incident fields  {\color{lxl}since} the fundamental solution $\Phi(x,y)=\frac{\i}{4}H_0^{(1)}(k|x-y|)$ belongs to the space $H_r^1(\Omega^p_H)$ only if $r<0$. From \cite{Lechl2017}, the direct solver introduced in Section 2.2 does not converge.
\end{remark}

The incident field $u^i\in H_r^1(\Omega^p_H)$ for any $r\in(0,1)$. 
Let $L=4$, then we use {\color{lxl}two incident fields $u^i$ and $u^i_L:=u^i(\cdot+2\pi L,\cdot)$}. Let $u^s$ and $u^s_L$ be the scattered fields {\color{lxl}corresponding} to the incident fields $u^i$ and $u^i_L$, and $u$, $u_L$ be the corresponding  {\color{lxl}total fields}.  From the estimation in Section 3, the error between $u_T^L:=u\circ\Phi_p$ and $u_0^L$, which is the total field with incident field $u^i_L$ and the periodic surface,  is bounded by:
\begin{equation*}
 \left\|u_T^L-u_0^L\right\|_{H^1(\Omega_H)}\leq C|8\pi|^{-r-1/2}\leq 0.008 C.
\end{equation*}
Note that the noise level is $\sigma =5\%$, $u_T^L$ could be treated as a good approximation of $u_0^L$ when the constant $C$ is assumed to be not too large.

\begin{figure}[H]
\centering
 \includegraphics[width=\textwidth]{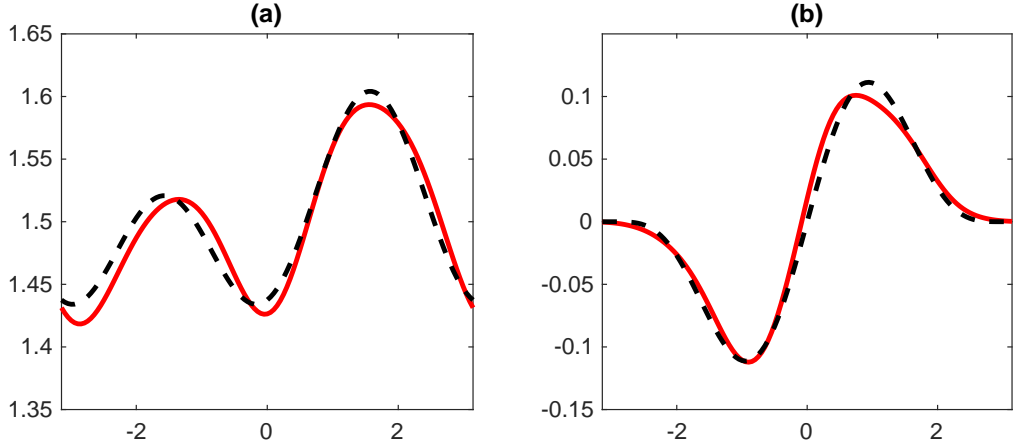}
 \caption{Example 1. (a): reconstruction of $\zeta_1$; (b): reconstruction of $p_1$. Black dotted curves: exact values; red curves: reconstructions.}
 \label{fig:eg1}
\end{figure}

\begin{figure}[H]
\centering
 \includegraphics[width=\textwidth]{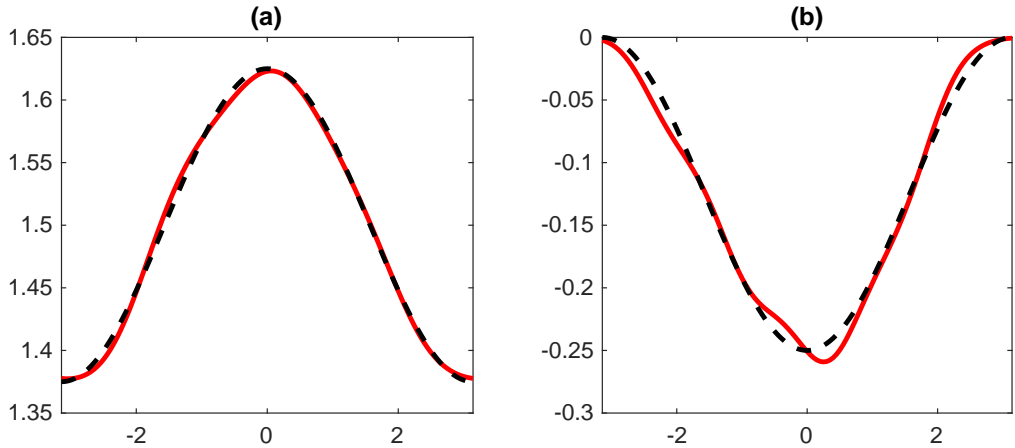}
 \caption{Example 2. (a): reconstruction of $\zeta_2$; (b): reconstruction of $p_2$. Black dotted curves: exact values; red curves: reconstructions.}
 \label{fig:eg2}
\end{figure}

Then we apply Algorithm \ref{alg} to reconstruct the perturbation  {\color{lxl}with} the known values $J$ and ${\bm c_0^1}$ from the sampling method. The reconstructs for Example 1  {\color{lxl}and Example 2} are shown in  {\color{lxl}Figure \ref{fig:eg1} and \ref{fig:eg2}, respectively}. From the left pictures of the two figures, the periodic surfaces are well reconstructed; based on the results for the periodic surfaces, we can also reconstruct the local perturbations very well.

\section*{Acknowlegdments} {\em This paper is devoted to Professor Armin Lechleiter. We will never forget him as a talented mathematician, a supportive colleague, and a dear friend.}

\bibliographystyle{alpha}
\bibliography{ip-biblio} 

\end{document}